\documentclass{article}
\usepackage{amsmath}
\usepackage{amsfonts,amsthm}
\usepackage{amssymb,enumerate}
\usepackage[pdftex]{graphicx}

\title{Almost all friendly matrices have many obstructions}
\date{}
\author{Richard Montgomery\footnote{Department of Pure Mathematics and Mathematical Statistics, Centre for Mathematical Sciences, Wilberforce Road, Cambridge, CB3 0WB, UK. r.h.montgomery@dpmms.cam.ac.uk}}

\newcounter{myfootnote}[page]

\newcounter{Lemma}
\newcounter{Theorem}
\newcounter{Cm}[section]

\newtheorem{lemma}[Lemma]{Lemma}

\newtheorem{corollary}[Lemma]{Corollary}
\newtheorem{theorem}[Theorem]{Theorem}

\theoremstyle{definition}
\newtheorem*{defn}{Definition}

\newtheorem{claim}[Cm]{Claim}

\theoremstyle{remark}
\newtheorem*{rmk}{Remark}

\newtheorem*{example}{Example}
\newtheorem*{nte}{Note}

\global\long\def\a{\alpha}
\global\long\def\b{\beta}

\global\long\def\t{\tau}
\global\long\def\s{\sigma}
\global\long\def\e{\varepsilon}

\global\long\def\N{\mathbb{N}}

\global\long\def\P{\mathbb{P}}

\global\long\def\G{\mathcal{G}}
\global\long\def\GG{\mathcal{G}}

\global\long\def\E{\mathbb{E}}

\global\long\def\TT{\mathcal{T}}

\global\long\def\PP{\mathcal{P}}

\global\long\def\re{\begin{rmk}}
\global\long\def\mark{\end{rmk}}
\global\long\def\ex{\begin{example}}
\global\long\def\ple{\end{example}}
\global\long\def\no{\begin{nte}}
\global\long\def\ted{\end{nte}}
\global\long\def\en{\begin{compactenum}}
\global\long\def\um{\end{compactenum}}
\global\long\def\li{\begin{compactitem}}
\global\long\def\st{\end{compactitem}}
\global\long\def\de{\begin{defn}}
\global\long\def\fn{\end{defn}}
\global\long\def\cor{\begin{corollary}}
\global\long\def\ary{\end{corollary}}
\global\long\def\lem{\begin{lemma}}
\global\long\def\ma{\end{lemma}}
\global\long\def\arr{\begin{array}}
\global\long\def\ay{\end{array}}
\global\long\def\pr{\begin{proof}}
\global\long\def\oof{\end{proof}}

\begin{document}
\maketitle

\begin{abstract}
A symmetric $m\times m$ matrix $M$ with entries taken from $\{0,1,\ast\}$ gives rise to a graph partition problem, asking whether a graph can be partitioned into $m$ vertex sets matched to the rows (and corresponding columns) of $M$ such that, if $M_{ij}=1$, then any two vertices between the corresponding vertex sets are joined by an edge, and, if $M_{ij}=0$, then any two vertices between the corresponding vertex sets are not joined by an edge. The entry $\ast$ places no restriction on the edges between the corresponding sets. This problem generalises graph colouring and graph homomorphism problems.

A graph with no $M$-partition but such that every proper subgraph does have an $M$-partition is called a minimal obstruction. Feder, Hell and Xie~\cite{FHX07} have defined friendly matrices and shown that non-friendly matrices have infinitely many minimal obstructions. They showed through examples that friendly matrices can have finitely or infinitely many minimal obstructions and gave an example of a friendly matrix with an NP-complete partition problem. Here we show that almost all friendly matrices have infinitely many minimal obstructions and an NP-complete partition problem.
\end{abstract}

\section{Introduction}
Many graph partition problems can be described by the following framework. Given a symmetric $m$ by $m$ matrix $M$ with entries taken from $\{0,1,\ast\}$, a graph $G$ has an $M$-partition if its vertex set can be partitioned into $m$ vertex classes $C_1,\ldots,C_m$ so that, if $M_{ij}=0$, then there are no edges between vertices in $C_i$ and vertices in $C_j$  and, if $M_{ij}=1$, then every edge is present between vertices in $C_i$ and vertices in $C_j$. The symbol $\ast$ places no restriction on edges between the two corresponding classes. The case $i=j$ is included, so that if, for example, $M_{ii}=1$, then $C_i$ must induce a complete graph.

This framework generalises graph colouring and homomorphism problems. Indeed, if an $m$ by $m$ matrix $M_m$ has diagonal entries $0$ and all off-diagonal entries $\ast$ then a graph has a $M_m$-partition exactly when it is $m$-colourable. Suppose the matrix $M_H$ is formed from a graph $H$ by matching a row and column to each vertex and letting an entry in a row and column be $\ast$ if the matching vertices are connected by an edge, and $0$ otherwise. Then a graph $G$ is $M_H$ partitionable exactly when there is a graph homomorphism from~$G$ into~$H$.

For a fixed matrix $M$, the computational problem of determining whether a graph is $M$-partitionable was introduced and studied by Feder, Hell, Klein and Motwani~\cite{FHKM}. It is well known that determining whether a graph is $k$-colourable is a polynomial problem if $k\leq 2$, and is NP-complete otherwise. Determining whether a graph has a homomorphism into a fixed graph $H$ is a polynomial problem if $H$ is bipartite, and is NP-complete otherwise~\cite{HN90}. It is unknown whether every $M$-partition problem is either polynomial or NP-complete~\cite{FHKNP, hell14}.

If a matrix $M$ has a $\ast$ on the diagonal then trivially all graphs have an $M$-partition and further constraints must be included (for example insisting each class is non-empty) to form an interesting problem. Here we will assume that $M$ has no $\ast$ on the diagonal.

The class of graphs which are $M$-partitionable, for some fixed matrix $M$, forms a hereditary property, that is a graph property closed under removing vertices. An equivalent notation for these properties was introduced by Bollob\'as and Thomason~\cite{BT99} who defined a type to be a complete graph with vertices coloured blue or red, and edges coloured red, blue or green. An embedding of a graph into a type embeds a complete, or empty, subgraph into each blue, or red, vertex respectively and a complete, or empty, bipartite graph across each blue, or red, edge respectively. Any collection of edges and non-edges can be mapped across a green edge. This notion of types is equivalent to using symmetric matrices with entries from $\{0,1,\ast\}$ with no $\ast$ on the diagonal, where a type has a vertex for each pair of corresponding rows and columns and the symbols $0,1$ and $\ast$ correspond to the colours red, blue and green respectively.

Work by Pr{\"o}mel and Steger~\cite{psIII}, Bollob\'as and Thomason~\cite{BT99}, and Thomason and Marchant~\cite{MTorig} built to the conclusion that any hereditary property could be approximated by the property of being $M$-embeddable, for some symmetric matrix $M$. That is, for any hereditary property $\PP$ and fixed probability $p$ there is some symmetric matrix $M$ for which
\[
\log( \P(\G_{n,p}\in \PP))=\log( \P(\G_{n,p} \text{ has an $M$ partition}))+o(1).
\]

A simple set of matrices with a polynomial partition problem are those for which there is some finite collection of graphs such that a graph has an $M$-partition exactly when it excludes each of those finitely many graphs as an induced subgraph, that is, a subgraph formed by deleting vertices.
In this case, we may check whether $G$ contains any of these finitely many graphs as induced subgraphs in polynomial time. 
As $M$-partitionable graphs form a hereditary property, the set of $M$-partitionable graphs may always be determined by a (possibly infinite) set of forbidden induced subgraphs. Minimal obstructions to an $M$-partition problem are graphs which have no $M$-partition yet any proper induced subgraph does have an $M$-partition. Any set of forbidden induced subgraphs defining the $M$-partition problem must contain all the minimal obstructions to an $M$-partition, and the minimal obstructions are sufficient to define the problem.

Feder, Hell and Xie~\cite{FHX07} used a random construction to demonstrate that there are infinitely many minimal obstructions for the $M$-partition problem if we can find either of the following submatrices in $M$ by taking two rows and their corresponding columns.
\[
\left(
\begin{array}{cc}
0 & \ast \\
\ast & 0 \\
\end{array}
\right)
\hspace{3cm}
\left(
\begin{array}{cc}
1 & \ast \\
\ast & 1 \\
\end{array}
\right)
\]
They defined `friendly' matrices as those matrices which do not have this submatrix property, so that we know `unfriendly' matrices have infinitely many minimal obstructions. They demonstrated the existence of friendly matrices which have finitely many minimal obstructions, those with infinitely many minimal obstructions yet a polynomial partition problem, and those which have an NP-complete partition problem. Feder, Hell and Shklarsky~\cite{FHS13} have demonstrated that any matrix only has finitely many minimal obstructions which are split graphs, that is graphs which can be partitioned into a clique and an empty set with no restrictions on the edges between them.

Given a friendly matrix $M$ it appears difficult to determine whether or not it has finitely many minimal obstructions. Neither the matrices with infinitely many minimal obstructions, nor those with finitely many minimal obstructions, form a class of matrices closed under deleting pairs of matching rows and columns. Let $M_i$ be the matrix formed from $M$ by deleting the $i$th row and the $i$th column and suppose $M$ is written as
\[
\left(\begin{array}{cc}
A & C \\
C & B \\
\end{array}
\right),
\]
where $A$ has each diagonal entry 0 and $B$ has each diagonal entry 1. Feder, Hell and Xie~\cite{FHX07} showed that if either $A$ or $B$ has no two rows the same and $M_i$ has finitely many minimal obstructions for each $i$ then $M$ itself has finitely many obstructions. If we select a friendly matrix randomly and uniformly from all $2n$ by $2n$ friendly matrices with $n$ entries for both 1 and 0 on the diagonal then almost surely this first condition on $A$ and $B$ holds. That is, as $n$ increases the probability this condition holds tends to $1$. However, we will show here that almost all friendly matrices have infinitely many minimal obstructions.

\begin{theorem}\label{aafriendly}
Almost all friendly matrices have infinitely many minimal obstructions.
\end{theorem}

Using the same notation, if $A$ and $B$ both have no three rows the same and $M_i$ has a polynomial partition problem for each $i$ then $M$ itself has a polynomial partition problem~\cite{FHX07}. While it appears to be difficult to determine whether or not a friendly matrix $M$ has an NP-complete partition problem, we will show that almost all friendly matrices do have an NP-complete partition problem.

\begin{theorem}\label{aafriendlyNP}
Almost all friendly matrices have an NP-complete partition problem.
\end{theorem}

While unfriendly matrices are known to have infinitely many minimal obstructions, examples occur both where the partition problem is polynomial and where it is NP-complete. Indeed, the two unfriendly matrices encoding the 2-colouring and 3-colouring problem respectively provide such examples. We will show that almost all matrices have an NP-complete partition problem, where we select a matrix uniformly at random from the $n$ by $n$ symmetric matrices with entries in $\{0,1,\ast\}$ which do not have a $\ast$ on the diagonal.

\begin{theorem}\label{aaNP}
Almost all matrices have an NP-complete partition problem.
\end{theorem}

\section{Definitions and notation}
As the work here has a probabilistic flavour it will be convenient to use the notation of types introduced by Bollob\'as and Thomason~\cite{BT99} in their work on hereditary graph properties.
\de
A \emph{type} $\tau$ is a complete graph where each vertex is coloured either red or blue and each edge is coloured red, blue or green. We will denote the set of red vertices by $R(\tau)$, the set of blue vertices by $B(\tau)$ and the entire vertex set by $V(\tau)$.
\fn
\de
An \emph{embedding} of a graph $G$ into a type $\tau$ is a map $\psi:V(G)\to V(\tau)$ where if $uv$ is an edge in $G$ then either $u$ and $v$ are mapped by $\psi$ to the same blue vertex, or $\psi(u)\psi(v)$ is a blue or green edge, and if $uv$ is not an edge in $G$ then either $u$ and $v$ are mapped by $\psi$ to the same red vertex, or $\psi(u)\psi(v)$ is a red or green edge.

If a graph $G$ has an embedding into the type $\tau$ then we say that it is \emph{embeddable into $\tau$}.
\fn
Thus, the vertices of the graph embedded into a red vertex form an independent set, and those embedded into a blue vertex form a clique. Only edges may be embedded across a blue edge. Only non-edges may be embedded across a red edge. We place no restriction on the edges and non-edges embedded across the green edges. 

For example, the graphs which are embeddable into the type with $k$ red vertices with green edges between them are exactly the $k$-colourable graphs. The $M$-partition problem is equivalent to a $\tau$-embedding problem for a type with a red vertex for each $i$ with $M_{ii}=0$ and a blue vertex for each $i$ with $M_{ii}=1$, with edges between vertices coloured red, blue or green if the corresponding entry between the rows and columns of $M$ is $0,1$ or $\ast$ respectively.

A \emph{friendly type} is defined analogously to a friendly matrix. That is, it is a type with no green edge between any two red vertices or between any two blue vertices.
\de
Given two types $\s$ and $\t$, an \emph{edge-homomorphism} $\phi:\s\to\t$ is a mapping of the vertices of $\s$ to the vertices of $\t$ where
\begin{itemize}
\item if $vw$ is a red edge, then either $\phi(v)=\phi(w)$ and this is a red vertex, or $\phi(v)\phi(w)$ is a red or green edge, and
\item if $vw$ is a blue edge, then either $\phi(v)=\phi(w)$ and this is a blue vertex, or $\phi(v)\phi(w)$ is a blue or green edge.
\end{itemize}
\fn
If in addition $\phi:\s\to\t$ preserves vertex colour and maps green edges across green edges then we say $\phi$ is a \emph{type-homomorphism}. If $G$ has an embedding $\psi$ into $\s$ and $\phi:\s\to\t$ is a type-homomorphism then $\phi\psi$ is an embedding of $G$ into $\t$. However, it will be convenient later to use the weaker notion of edge-homomorphism.

Finally, we must define a random friendly type and a random type.
\de
The random friendly type $\TT_f(n)$ is a type with $n$ red vertices and $n$ blue vertices with each edge between a red vertex and a blue vertex coloured red, green or blue uniformly and independently and each other edge coloured red or blue uniformly and independently.
\fn
\de
The random type $\TT(n)$ is a type with $n$ vertices where each edge is coloured red, green or blue uniformly and independently and each vertex is coloured red or blue uniformly and independently.
\fn
The first definition above is equivalent to choosing a matrix uniformly from the friendly $2n$ by $2n$ matrices with $n$ entries for both 1 and 0 on the diagonal. In our proofs of Theorems \ref{aafriendly} and \ref{aafriendlyNP} we will consider therefore the random friendly type $\TT_f(n)$. The second definition is equivalent to choosing a matrix uniformly from the symmetric $n$ by $n$ matrices with entries in $\{0,1,\ast\}$ and without a $\ast$ on the diagonal. In our proof of Theorem \ref{aaNP} we will consider therefore the random type $\TT(n)$.

\section{Properties of almost all random types}
For the proofs of the main theorems we will require some properties of almost all random types, which will be stated and proved here. A type $\s$ is a \emph{subtype} of the type $\tau$, denoted by $\s\subset\tau$, if $\s$ can be formed by deleting vertices from $\tau$. Where $A\subset V(\tau)$, we denote the subtype of $\tau$ with the vertex set $A$ by $\tau|_A$.

Given a subtype $\s\subset\t$ and a vertex $v\in V(\s)$, we say $v$ is a \emph{fixed point} of the function $\phi:V(\s)\to V(\t)$ if $\phi(v)=v$. As usual, we say a property of $\TT_f(n)$ or $\TT(n)$ holds \emph{with high probability} if it fails with probability $o(1)$ as $n\to\infty$.

\lem \label{typehomo}
Let $\a>\b>0$ and $\t=\TT_f(n)$ or $\TT(n)$. With high probability the following is true. For every subtype $\s\subset\tau$ with at least $\a n$ vertices and every edge-homomorphism $\phi:\s\to\t$, $\phi$ has at least $\b n$ fixed points.
\ma
\pr Let $\tau=\TT(n)$, using the vertex set $[n]=\{1,\ldots,n\}$. The case for $\TT_f(n)$ follows similarly. Let $\Phi$ be the set of pairs $(A,\phi)$ of vertex sets $A\subset [n]$ and functions $\phi:A\to [n]$ where $|A|\geq \a n$ and $\phi$ has at most $\beta n$ fixed points. As there are certainly at most $n^n$ such functions $\phi$ for each set $A$, we have $|\Phi|\leq 2^nn^n$. Let $X$ be the number of pairs $(A,\phi)\in \Phi$ for which $\phi$ is an edge-homomorphism from $\tau|_A$ to $\tau$. Thus, $X$ is a random variable dependant on $\tau$. To prove the lemma we need to show that $\P(X>0)\to 0$ as $n\to\infty$.

Let $(A,\phi)\in \Phi$. We will find an upper-bound for the probability that $\phi$ is an edge-homomorphism. We will first find a subset $A'\subset A$ with $|A'|\geq (\a-\b)n/4$ so that the image of $A'$ under $\phi$, $\phi(A')$, is disjoint from $A'$. Calculating the probability that $\phi|_{A'}$ is an edge-homomorphism will give us an upper-bound for the probability that $\phi$ is an edge-homomorphism.

Let $A_0=A$. Delete from $A_0$ any vertex fixed by $\phi$. As $(A,\phi)\in\Phi$, at least $(\a-\b)n$ vertices remain in $A_0$. Next, if there is a vertex $v$ for which there is an integer $k>1$ such that $\{v,\phi(v),\ldots,\phi^{k-1}(v)\}\subset A_0$ and $\phi^k(v)=v$, then delete the vertex $v$ from $A_0$. Repeat this process until no such vertex in $A_0$ exists. Note that after a vertex $v$ is deleted from $A_0$ all of the vertices $\phi(v),\ldots,\phi^{k-1}(v)$ will not subsequently be deleted from $A_0$. Therefore, at the end of this process, $A_0$ will still contain at least $(\a-\b)n/2$ vertices.

Take a new directed graph $H$ with the vertex set $A_0$ and edges $\vec{uv}$ exactly when $\phi(u)=v$. If, with the edge directions forgotten, there is a cycle in this graph, then it must be a directed cycle as each vertex has out-degree at most~1. If there is a cycle of length $k$ in $H$, then, taking some vertex $v$ in the cycle, we have that $\{v,\phi(v),\ldots,\phi^{k-1}(v)\}\subset A_0$ and $\phi^k(v)=v$. Therefore, as no such vertices remain in $A_0$, $H$ contains no cycles. The graph $H$ is therefore a forest, and hence is bipartite. As $|A_0|\geq (\a-\b)/2$, we can take $A'\subset A_0$ with $|A'|\geq\e n$, where $\e=(\a-\b)/4>0$, so that $A'$ is an independent set in $H$. Observe that, as $A'$ is an independent set, $\phi(A')$ and $A'$ are disjoint.

The probability that $\phi$ is an edge-homomorphism is at most the probability that $\phi|_{A'}$ is an edge-homomorphism. Let $k=|\phi(A')|$ and $\phi(A')=\{v_1,\ldots,v_k\}$, and, for each $i$, let $g_i$ be the number of vertices from $A'$ mapped to $v_i$ by $\phi$.

If $\phi$ is an edge-homomorphism, then if a vertex $v_i$ is coloured red (respectively, blue) it must have no blue (respectively, red) edges mapped into it, which has probability at most $\left(\frac23\right)^{\binom{g_i}{2}}\leq \left(\frac79\right)^{\binom{g_i}{2}}$. The probability an edge $v_iv_j$ is not green is $\frac{2}{3}$ whereupon the $g_ig_j$ edges mapped into it either cannot be red (if $v_iv_j$ is blue) or cannot be blue (if $v_iv_j$ is red), which certainly has probability at most $\frac23$ (as $g_ig_j\geq 1$). Therefore, for each of the $\binom{k}{2}$ edges $v_iv_j$ individually, the probability that $v_iv_j$ and the edges mapped into it are coloured so as not to prevent $\phi$ being an edge-homomorphism is at most $\frac13+\frac{2}{3}\times\frac23=\frac79$. Thus,
\[
\P(\phi\text{ is an edge-homomorphism})\leq \left(\frac79\right)^{\binom{k}{2}+\sum_{i=1}^k\binom{g_i}{2}}
\leq \left(\frac79\right)^{\binom{k}{2}+k\binom{\frac1k\sum_{i=1}^kg_i}{2}}
\]
\[
\leq \left(\frac79\right)^{\binom{k}{2}+k\binom{\frac{\e n}{k}}{2}}\leq\left(\frac79\right)^{C(k^2+\frac1kn^2)}\leq \left(\frac79\right)^{Cn^{4/3}}
\]
for some constant $C>0$, where in the last inequality we have considered the extremal case $k=2^{-1/3}n^{2/3}$.

As $|\Phi|\leq 2^nn^n$, we have, using Markov's inequality,
\[
\P(X>0)\leq \E(X)\leq 2^{n+n\log_2 n}\left(\frac79\right)^{C n^{4/3}}\to 0\text{\; as \;}n\to\infty.\qedhere
\]
\oof

In fact this proof is sufficient to show there are at most $C(n\log n)^{3/4}$ non-fixed vertices, for some constant $C>0$, but the statement of Lemma \ref{typehomo} is strong enough for our purposes. Furthermore, generalising the techniques used in proving Theorem \ref{aafriendly} would show that there can be at most $C\log n$ non-fixed vertices, for some constant $C>0$, which is tight up to the constant.

The techniques in the proofs of Theorems \ref{aafriendly}, \ref{aafriendlyNP} and \ref{aaNP} will require some technical lemmas, for which we need to define the common neighbourhood of a vertex set of a type.

\de
For a vertex subset $A$ of a type $\tau$,
\[
N(A)=\{v\in V(\tau)\setminus A\,:\text{ there do not exist }r,b\in A\text{ with }vr\text{ red and }vb\text{ blue}\}.
\]
When $A$ contains few vertices we list them, for example writing $N(v,w)$ for $N(\{v,w\})$.
\fn

\lem \label{nsize}
Given $\tau =\TT_f(n)$, with high probability it is true that for all distinct $r_1,r_2\in R(\tau)$ and distinct $b_1,b_2\in B(\tau)$ the following holds.
\begin{enumerate}[i)]
\item $|N(r_1,r_2)\cap N(b_1,b_2)|\geq \frac23 n$, and,
\item $|N(r_1,r_2)\cap N(b_1,b_2)\cap N(v,w)|\leq \frac{16}{27} n$ for all distinct $v,w\in V(\tau)$ with $\{v,w\}\neq\{r_1,r_2\}$ and $\{v,w\}\neq\{b_1,b_2\}$.
\end{enumerate}
\ma

\pr i) A red vertex $v\in R(\t)\setminus\{r_1,r_2\}$ is in $N(r_1,r_2)$ unless $vr_1$ and $vr_2$ take different colours in $\{\text{red},\text{blue}\}$, which, as $\tau$ is a random friendly matrix, happens with probability $\frac12$. A red vertex $v\in R(\t)\setminus\{r_1,r_2\}$ is in $N(b_1,b_2)$ unless $vb_1$ and $vb_2$ take different colours in $\{\text{red},\text{blue}\}$. There are thus 2 bad colourings out of the 9 different possible colourings of $vr_1$ and $vr_2$.
Therefore,
\[
\P(v\in N(r_1,r_2)\cap N(b_1,b_2)) =\frac12\left(1-\frac29\right)=\frac{7}{18}.
\]
Hence, $|N(r_1,r_2) \cap N(b_1,b_2)\cap R(\tau)|$ is a binomial variable with parameters $n-2$ and $\frac{7}{18}$. By Chernoff's inequality, we have for any $\e$, $0<\e\leq 1$, that
\[
\P\left(|N(r_1,r_2) \cap N(b_1,b_2)\cap R(\tau)|< \frac{7}{18}\left(1-\e\right)(n-2)\right)\leq \exp\left(-\frac{7}{72}\e^2(n-2)\right).
\]
By taking $\e=\frac1{8}$, we have, for $n\geq 100$,
\[
\P\left(|N(r_1,r_2) \cap N(b_1,b_2)\cap R(\tau)|< \frac{1}{3}n\right)\leq \exp\left(-\frac{1}{1000}n\right).
\]
The same result holds similarly for $|N(r_1,r_2) \cap N(b_1,b_2)\cap B(\tau)|$, and thus, for $n\geq 100$,
\[
\P\left(|N(r_1,r_2) \cap N(b_1,b_2)|< \frac{2}{3}n\right)\leq 2 \exp\left(-\frac{1}{1000}n\right).
\]
There are at most $n^4$ choices for the vertices $r_1,r_2,b_1,b_2$, so the expected number of sets of such vertices for which i) does not hold is at most $2n^4\exp(-\frac{1}{1000}n)$, if $n\geq 100$. Thus, with high probability all sets of such vertices will satisfy i).

ii) Proving the second part is very similar to proving the first part, but requires several different cases based on whether each of $v$ and $w$ is red or blue and whether each of $v$ and $w$ is in the set $\{r_1,r_2,b_1,b_2\}$ or not. In each case $|N(r_1,r_2)\cap N(b_1,b_2)\cap N(v,w)|$ can be shown to be a sum of two binomial variables whose combined expectation is strictly less than $\frac{16}{27}n$, for large $n$. Note that if a vertex lies in $N(r_1,r_2)\cap N(b_1,b_2)\cap N(v,w)$ it satisfies additional restrictions than if it just lies in $N(r_1,r_2)\cap N(b_1,b_2)$, and thus we expect the combined expectation in each case to be distinctly smaller than the corresponding combined expectation in part i).

The highest combined expectation is when $\{v,w\}\subset\{r_1,r_2,b_1,b_2\}$, or, without loss of generality, when $v=r_1$ and $w=b_1$. In this case, a red vertex $x\notin \{r_1,r_2\}$ is in $N(r_1,r_2)\cap N(b_1,b_2)\cap N(r_1,b_1)$ if $xr_1$ and $xr_2$ are the same colour and if the following holds. Either $xb_1$ is green, or $xb_1$ and $xb_2$ are both the same colour as $xr_1$, or $xb_1$ is the same colour as $xr_1$ and $xb_2$ is green. Thus, a red vertex not in $\{r_1,r_2\}$ is in $N(r_1,r_2)\cap N(b_1,b_2)\cap N(r_1,b_1)$ with probability $\frac12\left(\frac13+\frac19+\frac19\right)=\frac{5}{18}$. The probability is the same for a blue vertex not in $\{b_1,b_2\}$, so the sum of the expectations of the distributions is $\frac{5}{9}(n-2)$, strictly less than $\frac{16}{27}n$, for large $n$, as required.

In each case, Chernoff's inequality shows the probability that $|N(r_1,r_2)\cap N(b_1,b_2)\cap N(v,w)|>\frac{16}{27}n$ is exponentially small.
As there are at most $n^6$ choices for the vertices $r_1,r_2,b_1,b_2,v,w$, with high probability $\t$ will satisfy the second condition of the lemma.
\oof
The following two lemmas will take the same role in the proofs of Theorems \ref{aafriendlyNP} and \ref{aaNP} respectively as Lemma \ref{nsize} will take in the proof of Theorem \ref{aafriendly}. They can be proved similarly to Lemma \ref{nsize}.
\lem \label{nsize2}
Given $\tau =\TT_f(n)$, with high probability it is true that for every set $A\subset V(\tau)$ of six red and three blue vertices the following holds.
\begin{enumerate}[i)]
\item $|N(A)|\geq\frac{1}{36}n$, and
\item $|N(A\cup\{v\})|\leq\frac{1}{40}n$ for each vertex $v\notin A$.\hfill\qed
\end{enumerate}
\ma
\lem \label{nsize3}
Given $\tau =\TT(n)$, with high probability it is true that for every set $A\subset V(\tau)$ of three vertices the following holds.
\begin{enumerate}[i)]
\item $|N(A)|\geq\frac{14}{27}n$, and
\item $|N(A\cup\{v\})|\leq\frac{13}{27}n$ for each vertex $v\notin A$.\hfill\qed
\end{enumerate}
\ma

\section{Proof of the main results}
\pr[Proof of Theorem \ref{aafriendly}]
Let $\rho$ be the type with the vertex set $\{r_1,r_2,r_3,b_1,b_2,b_3\}$ so that
\begin{itemize}
\item the vertices $r_1,r_2$ and $r_3$ are red, and the vertices $b_1,b_2$ and $b_3$ are blue,
\item the edges $r_1r_3$, $r_2r_3$ and $b_1b_2$ are blue,
\item the edges $r_1b_1$, $r_1b_3$, $r_2b_2$ and $r_3b_2$ are green, and
\item all the remaining edges within $\rho$ are red.
\end{itemize}
This subtype $\rho$ is the example used by Feder, Hell and Xie~\cite{FHX07} to demonstrate that there are friendly types with infinitely many minimal obstructions.

Let $\tau=\TT_f(n)$. Taking any set of three red and three blue vertices in $V(\tau)$, the probability that restricting $\tau$ to this set gives a copy of $\rho$ is at least $(\frac13)^{\binom 62}$. Considering $\lfloor \frac{n}{3}\rfloor$ disjoint such vertex sets in $V(\tau)$, and using Chernoff's inequality, demonstrates that $\tau$ will almost surely contain a subtype which is a copy of $\rho$. Relabelling, let us suppose that $\rho\subset \tau$.
By Lemma \ref{typehomo} and Lemma \ref{nsize}, $\tau$ almost surely also satisfies the conclusion of Lemma \ref{typehomo} with $\a=\frac23$ and $\b=\frac{17}{27}$, and the conclusion of Lemma \ref{nsize}. Assume additionally that $n\geq 100$. We will show that, for each $m\in\N$, $\tau$ has a minimal obstruction with at least $m$ vertices. The type $\tau$ must then have infinitely many obstructions, and the theorem will be proved.

Let $\s\subset\tau$ have the vertex set $N(r_1,r_2)\cap N(b_1,b_2)$. Then $r_3,b_3\in V(\sigma)$ and, by the property from Lemma \ref{nsize}, $|\s|\geq\frac23n$.
For each vertex $v\in V(\s)$, create a new vertex $v'$. Let $G$ be the graph of order $|V(\s)|+2m$ with vertex set $\{v':v\in V(\s)\}\cup\{x_1,\ldots,x_m,y_1,\ldots,y_m\}$, whose edge set comprises 
\begin{itemize}
\item the edges $u'v'$ for which $uv$ is a blue edge of $\s$,
\item the edges in the complete graph on $\{y_1,\ldots,y_m\}$,
\item the edges in the path $x_1y_1x_2\ldots x_my_m$,
\item the edges $x_1b_3'$ and $y_mr_3'$,
\item the edges in $\{x_iv' : v\in V(\sigma),\text{ }i\in[m] \text{ and }r_1v\text{ or }r_2v\text{ is blue in }\tau\}$, and
\item the edges in $\{y_iv' : v\in V(\sigma),\text{ }i\in[m] \text{ and }b_1v\text{ or }b_2v\text{ is blue in }\tau\}$.
\end{itemize}
We claim the resulting graph $G$ has the following two properties.

\begin{claim} \label{minsb}
$G$ does not have a $\tau$-embedding.
\end{claim} 
\begin{claim} \label{path}
For each $i$, $1\leq i\leq m$, $G-x_i$ does have a $\tau$-embedding.
\end{claim}
Claim \ref{minsb} implies that $G$ contains a minimal obstruction. Claim \ref{path} implies that any minimal obstruction contained in $G$ must contain each vertex $x_i$. Thus, $\tau$ has a minimal obstruction with at least $m$ vertices, as required.

To prove Claim \ref{minsb}, suppose to the contrary that $\psi:V(G)\to V(\tau)$ is some embedding of the graph $G$. Let $\phi:\s\to \tau$ be defined by $\phi(v)=\psi(v')$. We know that $|V(\s)|\geq \frac23 n$, and can observe that $\phi$ is an edge-homomorphism. Indeed, if $vw$ is a red edge in $\sigma$ then $v'w'\notin E(G)$ and hence as $\psi$ is an embedding either $\phi(v)=\psi(v')=\psi(w')=\phi(w)$ is a red vertex or $\psi(v')\psi(w')=\phi(v)\phi(w)$ is a red or green edge. The same follows with `blue' and `$v'w'\in E(G)$' in place of `red' and `$v'w'\notin E(G)$' respectively.

Therefore, by the property from Lemma \ref{typehomo}, $\phi$ has at least $\frac{17}{27}n$ fixed vertices. Call this set of fixed vertices $F$. For $v\in V(\s)$ and $f\in F\setminus\{v,\phi(v)\}$, if $vf$ is blue then $\phi(v)\phi(f)=\phi(v)f$ cannot be red, and if $vf$ is red then $\phi(v)f$ cannot be blue. Therefore $f\in N(v,\phi(v))$, for each $f\in F\setminus\{v,\phi(v)\}$. As $F\subset N(r_1,r_2)\cap N(b_1,b_2)$, we have
\[
|N(r_1,r_2)\cap N(b_1,b_2)\cap N(v,\phi(v))|\geq \frac{17}{27}n-2.
\]
Therefore, by the property from Lemma \ref{nsize}, to avoid a contradiction we must have $v=\phi(v)$, for each $v\in V(\s)$ (as $v\notin\{r_1,r_2,b_1,b_2\}$), and so $\phi$ in fact fixes every vertex in $\s$.
In other words, each vertex $v'$ in $G$ has been embedded into $v$ by $\psi$.

We now consider where the remaining vertices in $G$ are, that is the vertices in the path $x_1y_1x_2\ldots x_my_m$.
Fixing $i$, let $x=\psi(x_i)$. If $v\in V(\sigma)\setminus\{x\}$ and $r_1v$ is a blue edge in $\tau$, then $x_iv'$ is an edge in $G$. As $\psi(v')=v$, $x\neq v$, and $\psi$ is an embedding of $G$, we must have that $xv$ is a blue or green edge in $\tau$. Similarly, if $v\in V(\sigma)\setminus\{x\}$ and $r_1v$ is a red edge in $\tau$, then $xv$ is a red or green edge in $\tau$. Therefore, for each $v\in V(\sigma)\setminus \{x\}$, $v\in N(r_1,x)$, and hence $|N(r_1,r_2)\cap N(b_1,b_2)\cap N(r_1,x)|\geq \frac{2}{3}n-1$. Using the property from Lemma \ref{nsize}, $x$ must be $r_1$ or $r_2$. Therefore, for each vertex $x_i$, we have $\psi(x_i)\in\{r_1,r_2\}$.

Similarly, for each $y_i$, we must have $V(\sigma)\setminus\{\psi(y_i)\}\subset N(b_1,\psi(y_i))$, so that $|N(r_1,r_2)\cap N(b_1,b_2)\cap N(b_1,\psi(y_i))|\geq \frac{2}{3}n-1$. By the property from Lemma \ref{nsize} then, $\psi(y_i)\in\{b_1,b_2\}$.


We know that $x_1b_3'\in E(G)$, $b_3'$ is embedded into $b_3$ (as $b_3'\in V(G)$) and $r_2b_3$ is a red edge in $\tau$. Therefore, $\psi(x_1)\neq r_2$, and hence $\psi(x_1)=r_1$. Then, as $x_1y_1\in E(G)$ and $r_1b_2$ is a red edge in $\tau$ we can deduce $\psi(y_1)=b_1$. Continuing like this, and considering $x_2$ next, the path $x_1y_1x_2\ldots x_my_m$ must be embedded between $r_1$ and $b_1$, until finally we deduce that $\psi(y_m)=b_1$. But then the edge $y_mr_3'$ of $G$ is embedded across the red edge $r_3b_1$ in $\tau$, a contradiction. Thus no such embedding $\psi$ exists and Claim \ref{minsb} has been proved.

To prove Claim \ref{path} we can define, for each $i$, an embedding $\psi:G-x_i\to\tau$ as follows. For $v\in V(\s)$, let $\psi(v)=v'$. For each $j$, $j< i$, let $\psi(x_j)=r_1$ and $\psi(y_j)=b_1$. For each $j>i$, let $\psi(x_j)=r_2$ and $\psi(y_j)=b_2$. Let $\psi(y_i)=b_2$.

As removing $x_i$ breaks the path $x_1y_1\ldots x_my_m$, the contradiction that arose in proving Claim \ref{minsb} does not occur, and $\psi$ can easily be seen to be an embedding from the definition of $G$.
\oof

\pr[Proof of Theorem \ref{aafriendlyNP}]
Suppose $\t$ is a friendly type with $n$ blue and $n$ red vertices, where $n\geq 1000$, for which the conclusion of Lemma~\ref{nsize2} holds, as well as the conclusion of Lemma~\ref{typehomo} with $\a=\frac{1}{36}$ and $\b=\frac{1}{38}$. Feder, Hell and Xie~\cite{FHX07} showed that there is a friendly type, which we shall call $\rho$, which has 6 red and 3 blue vertices and an NP-complete partition problem. Suppose further that $\t$ contains this type $\rho$. We will show that the problem of determining whether $G$ has an embedding into $\rho$ can be solved by determining whether there is an embedding into $\t$ of a specific graph $G'$ constructed from $G$ in polynomial time. Therefore, the problem of determining whether a graph has an embedding into $\tau$ must also be NP-complete. Using a similar argument to that used in the proof of Theorem~\ref{aafriendly}, we can almost surely assume that $\rho\subset\TT_f(n)$. The other properties of $\tau$ follow almost surely for $\TT_f(n)$ by Lemmas~\ref{typehomo} and \ref{nsize2}, and therefore the random friendly type $\TT_f(n)$ will almost surely have an NP-complete embedding problem.

Given a graph $G$, extend it to form the graph $G'$ as follows. Let $\sigma$ be the subtype of $\tau$ with the vertex set $N(V(\rho))$. For each $v\in V(\sigma)$ create a vertex $v'$ of $G'$ and connect it with an edge to every vertex in $G$ exactly when at least one of the edges from $v$ to $\rho$ in $\tau$ is blue. Add the edge $v'w'$ to $G'$ if $vw$ is a blue edge in $\tau$. Note that $|N(V(\rho))|\geq \frac{1}{36}n$ by the property from Lemma~\ref{nsize2}.

Given an embedding $\psi$ of $G'$ into $\tau$, the edge-homomorphism $\phi:\sigma\to\tau$ defined by $\phi(v)=\psi(v')$ must have at least $\frac{1}{38}n$ fixed vertices, by the property of $\tau$ from Lemma \ref{typehomo}. Let $F$ be the set of these fixed vertices, and let $v\in V(G)$. Given $f\in F\setminus\{\psi(v)\}$, if there is some vertex $x\in V(\rho)$ for which $fx$ is blue (so that all the edges from $V(\rho)$ to $f$ must be blue or green as $f\in V(\s)$), then $f'v$ is an edge in $G'$, so that $\psi(f')\psi(v)=f\psi(v)$ is a blue or green edge in $\tau$. All the edges from $V(\rho)\cup\{\psi(v)\}$ to $f$ are therefore blue or green, and hence $f\in N(V(\rho)\cup\{\psi(v)\})$. If there are no edges between $V(\rho)$ and $f$ which are blue, then $f'v$ is not an edge in $G$. Hence $f\psi(v)$ is a red or green edge in $\tau$, and $f\in N(V(\rho)\cup\{\psi(v)\})$. Therefore, for each $f\in F\setminus\{\psi(v)\}$, $f\in N(V(\rho)\cup\{\psi(v)\})$, and hence $F\setminus\{\psi(v)\}\subset N(V(\rho)\cup\{\psi(v)\})$. Therefore, $|N(V(\rho)\cup\{\psi(v)\})|\geq \frac{1}{38}n-1$, and so, for large $n$, we must have by the property from Lemma \ref{nsize2} that $\psi(v)\in V(\rho)$. As this holds for all $v\in V(G)$, $\psi|_{V(G)}$ is an embedding of $G$ into $\rho$. Therefore, if the graph $G'$ has an embedding into $\tau$ then $G$ must have an embedding into $\rho$.

Given an embedding of $G$ into $\rho$ we can extend this to an embedding of $G'$ into $\tau$ by mapping $v'$ into $v$ for each $v\in V(\s)$. Therefore, as required, $G$ has an embedding into $\rho$ if and only if $G'$ has an embedding into $\tau$.
\oof

\pr[Proof of Theorem \ref{aaNP}] Let $\rho$ be the type consisting of three red vertices with a green edge between each pair of these vertices. Given a random type $\tau=\TT(n)$ we have, almost surely, that $\rho\subset \tau$ (with relabelling) and that the conclusion of Lemma \ref{nsize3} and Lemma \ref{typehomo}, with $\a=\frac{14}{27}$ and $\b=\frac{1}{2}$, hold for $\tau$ (as shown by those lemmas). Given a graph $G$, we may create from $G$ an auxiliary graph $G'$ in polynomial time which is embeddable into $\tau$ precisely when $G$ is embeddedable into $\rho$, using an identical construction to that used in the proof of Theorem \ref{aafriendlyNP}. Therefore, asking whether $G'$ has a $\tau$-embedding is equivalent to asking whether $G$ has a 3-colouring. The 3-colouring problem is known to be NP-complete, and thus the $\tau$-embedding problem must be NP-complete as well.
\oof

\bibliographystyle{plain}
\bibliography{rhmreferences}

\end{document}